\documentclass[A4paper,11pt]{article}
\usepackage{latexsym}
\usepackage{mathrsfs}
\usepackage{amssymb}
\usepackage{amscd}
\usepackage[dvips]{graphicx}                  

\newtheorem{theorem}{Theorem}

\newtheorem{proposition}{Proposition}

\newenvironment{proof}[1][Proof]{\textbf{#1.} }{\ \rule{0.5em}{0.5em}}
\long\def\symbolfootnote[#1]#2{\begingroup%
	\def\thefootnote{$\;$}\footnote[#1]{$^*$#2}\endgroup}
\begin{document}
	
	\title{Decreasing chains without lower bounds in Rudin-Frol\'ik order for regulars}
	\author{Joanna Jureczko}
\maketitle

\symbolfootnote[2]{Mathematics Subject Classification: Primary 03E10, 03E20, 03E30..

\hspace{0.2cm}
}

\begin{abstract}
	The aim of this paper is to prove that for  there exist a chain in the Rudin-Frol\'ik order of $\beta\kappa\setminus \kappa$ of length $\mu$ with $\kappa \leqslant \mu \leqslant 2^\kappa$ for regular $\kappa> \omega$ without a lower bound.
\end{abstract}

\section{Introduction}

 The Rudin-Frol\'ik order, which has been defined by Z. Frol\'ik in \cite{ZF}, is an important topic, but still little known. Initially, this order was used by Frol\'ik   to prove that $\beta\omega\setminus \omega$ is not homogeneous.  M.E. Rudin, who nearly defined this ordering  in \cite{MR1},  as the first observed that the relation between filters she used is really  ordering.
This order was later investigated by D. Booth in \cite{DB} who showed that this relation is a partial ordering of the equivalence classes, that is a tree, and that it is not well-founded.

The papers which are included in the scope of our considerations around the Rudin-Frol\'ik order are E. Butkovi\v cov\'a who worked on this topic between 1981 and 1990 publishing a number of papers concerning ultrafilters  in Rudin-Frol\'ik order in $\beta \omega\setminus \omega$.  
Let us briefly review her achievements in this topic.
In \cite{BB} with L. Bukovsk\'y and in \cite{BE3} she constructed an ultrafilter on $\omega$ with the countable set of its predecessors. In \cite{BE1} she  constructed ultrafilters without immediate predecessors. In \cite{BE3}, Butkovi\v cov\'a showed that there exists in Rudin-Frol\'ik order an unbounded chain  orded-isomorphic in $\omega_1$.  In \cite{BE4}, she proved that  there is a set of $2^{2^{\aleph_0}}$ ultrafilters incomparable in Rudin-Frol\'ik order which is bounded from below and no its subset of cardinality more than one has an infimum. In \cite{BE5}, Butkovi\v cov\'a proved that for every cardinal between $\omega$ and $\mathfrak{c}$ there is a strictly decreasing chain without a lower bound.

In 1976 A. Kanamori published a paper \cite{AK} in which, among others showed that the Rudin-Frol\'ik tree cannot be very high if one consider it over a measurable cardinal. Moreover, in the same paper he left a number of open problems about Rudin-Frol\'ik order. 
Recently,  M. Gitik in \cite{MG} answered some of them but using metamathematical methods. The solution of some of the problems from \cite{AK} presented in combinatorial methods are in preparation, (\cite{JJ_kanamori}).

However, the Rudin-Frol\'ik order was investigated mainly for $\beta\omega$, significant results may be obtained when considering this order for the space $ \beta \kappa $, where $ \kappa $ is any cardinal. 
Since some of the results are based on the construction of the sequences of filterx  by transfiite induction, a different technique is needed in the case of $\beta \kappa$ when $\kappa > \omega$.

The method proposed in \cite{BK} by Baker and Kunen  comes in handy
In mentioned paper the authors presented  very usefull method which can be recognized as a generalization of method  presented in \cite{KK}. 
It is worth empahsizing that both methods, (from \cite{KK} and \cite{BK}), provide usefull "technology"  for keeping the transfinite construction for an ultrafilter not finished before $\mathfrak{c}$ steps, (see \cite{KK}), and $2^\kappa, $  for $\kappa$ being infinite cardinal, (see \cite{BK}), but the second method has some limitations, among others $\kappa$ must be regular.
Due to the lack of adequate useful method for a singular cardinal $\kappa$, the similar results but for singulars are still left as open questions. So far, we have not found an answer whether the assumptions can be omitted, which would probably also involve  changes in the methods used in our considerations. Therefore, based on the results from \cite{BK}, we restrict our results to this particular case.

The motivation for this paper follows from two papers \cite{VD} and \cite{BE6}. In \cite{VD} the author constructed a sequence $\langle X_\beta \colon \beta < \mathfrak{c}\rangle$ of countable discrete subsets of $\beta \omega$ such that $X_\eta \subseteq \overline{X_\varepsilon} \setminus D_\eta$ whenever $\varepsilon < \eta < \mathfrak{c}$ and $|\bigcap_{\beta< \mathfrak{c}} \overline{D_\beta}|=1.$ In \cite{BE6} there is shown that there exists such a set of length $\mu < \mathfrak{c}$. As was shown in \cite{BE6} this topic is very closed to the Rudin-Frol'ik order of ultrafilters on $\omega$.

On the other hand the generalization   a number of results concerning chains and sets of $\beta\kappa$ for regular $\kappa$ in the Rudin-Frol\'ik order is presented in \cite{JJ_kanamori, JJ_order, JJ_order1, JJ_order2}. It seems to be natural to ggeneralize also the above result. Thus, the main result of this paper is the following theorem
\\\\
\textbf{Theorem 1.}
\textit{Let $\kappa$ be a regular cardinal and let $\hat{\varphi}$ be a $\kappa$-shrinking function. For every cardinal $\mu$ with $\kappa \leqslant \mu \leqslant 2^\kappa$ there exists a strictly decreasing chain in the Rudin-Frol\'ik order $\{\mathcal{G}_\beta \colon \beta < \mu\}$ without a lower bound}.
\\\\

Let us note that in \cite{BB} the above theorem is proved for $\mu = \omega$ and in \cite{VD} it is proved for $\mu = \mathfrak{c}$. In \cite{BE6} the above result is proved for $\kappa = \omega$.
On the other hand, using slightly different argumentation there is shown in \cite{JJ_order1} that there exists a chain order isomorphic to $(2^\kappa)^+$ for regular $\kappa$.
We assume that $cf(\mu) > \omega$, because the case $cf(\mu) = \omega$  is proved in \cite{BE6}.

The paper is organized as follows: in Section 2, there are presented definitions and previous facts needed for the results presented in further parts of this paper. In Section 3, there are proved auxiliary results, the main result and the open problem.

We have tried to present all the necessary definitions, assuming tacitly that the reader has a basic knowledge of ultrafilters and the Rudin-Frol\'ik order.
However, for definitions and facts not quoted here, I refer the reader to e.g. \cite{TJ, CN}.

\section{Definitions and previous results}

\textbf{2.1.}
In the whole paper, we assume that $\kappa$ is an infinite cardinal. Then $\beta\kappa$ means the \v Cech-Stone compactification, where $\kappa$ has the discrete topology. Hence, $\beta\kappa$ is the space of ultrafilters on $\kappa$ and $\beta \kappa \setminus \kappa$ is the space of nonprincipal ultrafilters on $\kappa$. 
\\\\
\textbf{2.2.}
A set $\{\mathcal{F}_\alpha \colon \alpha < \kappa\}$ of filters on $\kappa$ is \textit{$\kappa$-discrete} iff there is a partition $\{A_\alpha \colon \alpha < \kappa\}$ of $\kappa$ such that $A_\alpha \in \mathcal{F}_\alpha$ for each $\alpha < \kappa$. 
\\\\
\textbf{2.3.}  Let $\mathcal{F}, \mathcal{G}\in \beta \kappa\setminus \kappa$. We define \textit{Rudin-Frol\'ik order} as follows
$$\mathcal{F} \leqslant_{RF} \mathcal{G} \textrm{ iff } \mathcal{G} = \Sigma(X, \mathcal{F})$$
for some $\kappa$-discrete set $X = \{\mathcal{F}_\alpha \colon \alpha < \kappa\} \subseteq \beta\kappa,$
where $$\Sigma(X, \mathcal{F}) = \{A \subseteq \kappa \colon \{\alpha < \kappa \colon A \in \mathcal{F}_\alpha\}\in \mathcal{F}\}.$$
We define $$\mathcal{F} =_{RF} \mathcal{G} \textrm{ iff } \mathcal{F} \leqslant_{RF} \mathcal{G} \textrm{ and } \mathcal{G} \leqslant_{RF} \mathcal{F}$$  $$\mathcal{F} <_{RF} \mathcal{G} \textrm{ iff } \mathcal{F} \leqslant_{RF} \mathcal{G} \textrm{ and } \mathcal{F} \not =_{RF} \mathcal{G}.$$
Conversely, if $\mathcal{G} \in \overline{X}$ then there exists a unique ultrafilter $\Omega(X, \mathcal{G})$ such that $\Sigma(X, \Omega(X, \mathcal{G})) = \mathcal{G}$.
\\\\
\textbf{2.4.} 	Let us accept the following notation: 
\begin{itemize}
	\item $\mathcal{FR}(\kappa) = \{A \subset \kappa \colon |\kappa\setminus A|< \kappa\},$
	\item $[A, B, C,...]$ means a filter generated by $A, B, C, ...$.
\end{itemize}
\textbf{2.5.} A function $\hat\varphi \colon [\kappa^+]^{<\omega} \to [\kappa]^{< \omega}$ is \textit{$\kappa$-shrinking} iff 
\begin{itemize}
	\item [(1)] $p\subseteq q$ implies $\hat{\varphi}(p) \subseteq \hat{\varphi}(q)$, for any $p, q \in [\kappa^+]^{<\omega}$,
	\item [(2)] $\hat{\varphi}(0) = 0$.
\end{itemize}

A \textit{step family} (over $\kappa$, with respect to $\hat{\varphi}$) is a family of subsets of $\kappa$, $$\{E_t \colon t \in [\kappa]^{<\omega}\} \cup \{A_\alpha \colon \alpha < \kappa^+\}$$ satisfying the following conditions:
\begin{itemize}
	\item [(1)] $E_s\cap E_t = \emptyset$ for all $s, t \in [\kappa]^{<\omega}$ with $s \not =t$,
	\item [(2)] $|\bigcap_{\alpha \in p} A_\alpha \cap  \bigcup_{t\not \supseteq \hat{\varphi}(p)}E_t| < \kappa$ for each $p \in [\kappa^+]^{<\omega}$,
	\item [(3)] if $\hat{\varphi}(p) \subseteq t$, then $|\bigcap_{\alpha\in p} A_\alpha \cap E_t|=\kappa$ for each $p \in [\kappa^+]^{<\omega}$ and $t \in [\kappa]^{<\omega}$. 
\end{itemize}

Let $I$ be an index set and $\mathcal{F}$ be a filter on $\kappa$. The family $$\{E_t^i \colon t \in [\kappa]^{<\omega}, i \in I\} \cup \{A_\alpha^i \colon \alpha < \kappa^+, i \in i\}$$ is an \textit{independent matrix of $|I|$ step-families} (over $\kappa$) with respect to $\mathcal{F}, \hat{\varphi}$ iff
\begin{itemize}
	\item [(1)] for each fixed $i \in I$, $\{E_t^i \colon t \in [\kappa]^{<\omega}\} \cup \{A_\alpha^i \colon \alpha < \kappa^+\}$ is a step-family,
	\item [(2)] if $n \in \omega, p_0, p_1, ..., p_{n-1} \in [\kappa^+]^{<\omega}, t_0, t_1, ..., t_{n-1} \in [\kappa]^{<\omega}$,  $i_0, i_1, ..., i_{n-1} \in I$ with $i_k\not = i_m, k\not = m$ and $\hat{\varphi}(p_k) \subseteq t_k$, then 
	$$\bigcap_{k=1}^{n-1}(\bigcap_{\alpha\in p_k}A^{i_k}_{\alpha} \cap E^{i_k}_{t_k}) \in \mathcal{F}^+,$$
	where $\mathcal{F}^+ = \{D\subseteq \kappa \colon \kappa\setminus D \not \in \mathcal{F}\}$.  
\end{itemize}
\noindent
\textbf{Fact 1 (\cite{BK}).} If $\kappa$ is a regular cardinal and $\hat{\varphi}$ is a $\kappa$-shrinking function, then there exists and independent matrix of $2^\kappa$ step-families over $\kappa$ with respect to the  filter $\mathcal{RF}(\kappa)$, $\hat{\varphi}$.

\section{Main result and open problem}

The main result of this paper is to prove the following theorem.

\begin{theorem}
		Let $\kappa$ be a regular cardinal and let $\hat{\varphi}$ be a $\kappa$-shrinking function. For every cardinal $\mu$ with $\kappa \leqslant \mu \leqslant 2^\kappa$ there exists a strictly decreasing chain in the Rudin-Frol\'ik order $\{\mathcal{G}_\beta \colon \beta < \mu\}$ without a lower bound.
\end{theorem}

Before the starting the proof of Theorem 1 we will make some remarks.
SInce the proof of Theorem 1 is slightly based on the proof presented in \cite{BE6} we proceed as follows. 
We will construct a sequence $\langle X_\beta \colon \beta < \mu\rangle$ of $\kappa-$discrete subsets an a point $\mathcal{F} \in \beta \kappa \setminus \kappa$ such that 
\begin{itemize}
	\item [(1)] $\mathcal{F} \in \bigcap_{\beta < \mu} \overline{X_\beta}$,
	\item [(2)] $X_\beta \subseteq \overline{X_\beta} \setminus X_\gamma$ whenever $\gamma < \beta \mu$,
	\item [(3)] for every $\kappa-$discrete $Y \subseteq \bigcap_{\beta < \mu} \overline{X_\beta}\setminus \{\mathcal{G}\}$ we have $\mathcal{G} \not \in \overline{Y}$. 
\end{itemize}
Then $\langle \Omega(X_\beta, \mathcal{F}) \colon \beta < \mu\rangle$ will be strictly decreasing in the Rudin-Frol\'ik order.

For proving our main result we still need the following result, compare \cite{BE6}. The proof of the following proposition can be carried out analogically to the proof of this result in the case of $\kappa = \omega$ presented in \cite[Lemma, p. 254]{BE6}.

\begin{proposition}
	The sequence $\langle X_\beta \colon \beta < \mu\rangle$ is strictly decreasing in the Rudin-Frol\'ik order iff for every $\kappa-$discrete set $Y \subseteq \bigcap_{\beta < \mu} \overline{X_\beta}\setminus \{\mathcal{G}\}$ we have $\mathcal{G} \not \in \overline{Y}$. 
	\end{proposition}

\begin{proof}\textbf{of Theorem 1}
	By Fact 1 fix a matrix of $2^\kappa$ step families (over $\kappa$)
		$$\{E_t^i \colon t \in [\kappa]^{<\omega}, i \in 2^\kappa\} \cup \{A_\eta^i \colon \eta < \kappa^+, i \in 2^\kappa\}$$
		which is independent with respect to  $\mathcal{FR}(\kappa), \hat{\varphi}$. 
		
			For our purpose, we slightly modify this matrix by "shrinking" $A^i_\alpha$ to $A^i_\alpha \subseteq \bigcup \{E_t^i \colon t \in [\kappa]^{<\omega}\}$ and then expanding $E^i_t$ so that $\{E^i_t \colon t \in [\kappa]^{< \omega}\}$ is a partition of $\kappa$. 
		Thus, we obtain the matrix fullfilling the following conditions:
		\begin{itemize}
			\item [(a)] $E^i_s\cap E^i_t = \emptyset$ for all $s, t \in [\kappa]^{<\omega}$ with $s \not =t$,
			\item [(b)] $|\bigcap_{\eta \in p} A^i_\eta \cap  \bigcup_{t\not \supseteq \hat{\varphi}(p)}E^i_t| < \kappa$ for each $p \in [\kappa^+]^{<\omega}$,
			\item [(c)] if $\hat{\varphi}(p) \subseteq t$, then $|\bigcap_{\eta\in p} A^i_\eta \cap E^i_t|=\kappa$ for each $p \in [\kappa^+]^{<\omega}$ and $t \in [\kappa]^{<\omega},$ 
			\item [(d)] $\bigcup\{E^i_t \colon t \in [\kappa]^{< \omega}\} = \kappa.$
		\end{itemize}
The condition $(b)$ is still preserved after expanding $\{E_t^i \colon t \in [\kappa]^{<\omega}\}$ to a partition of $\kappa$. 
	\\Indeed. If there are $p_0 \in [\kappa^+]^{<\omega}$ and $i_0 \in 2^\kappa$ such that
	$$|\bigcap_{\eta \in p_0} A_\eta^{i_0} \cap  \bigcup_{t\not \supseteq \hat{\varphi}(p_0)}E_t^{i_0}| = \kappa,$$
	then $|\bigcup_{t\not \supseteq \hat{\varphi}(p_0)}E_t^{i_0}| = \kappa.$
	Then, by $(d)$ and $(a)$, there would exist  $t_0 \supseteq \hat{\varphi}(p_0)$ such that  $|E_{t_0}^{i_0}| < \kappa$. Hence
	$$|\bigcap_{\eta\in p_0} A^{i_0}_\eta \cap E^{i_0}_{t_0}|<\kappa.$$
	which contradicts $(c)$. 
	
	Now, we will construct the ultrafilters $\mathcal{G}_{\mu, t}$ and $\mathcal{\xi, t }$for any $\xi < \mu$ and $t \in [\kappa]^{< \omega}$ which fulfill conditions $(1)-(3)$. For this purpose we will construct the increasing sequences 
	$$\{\mathcal{F}^\alpha_{\xi, t} \colon \xi< \mu, t \in [\kappa]^{<\omega}\}$$ and 
	$$\{\mathcal{F}^\alpha_{\mu, t} \colon \xi< \mu, t \in [\kappa]^{<\omega}\}.$$
	Then, we will take $\mathcal{G}_{\xi, t} = \bigcup_{\alpha < 2^\kappa} \mathcal{F}^\alpha_{\xi, t}$ and $\mathcal{G}_{\mu, t} = \bigcup_{\alpha< 2^\kappa} \mathcal{F}^\alpha_{\mu, t}$.
	Then for each $mu$ with $\kappa \leqslant \mu \leqslant 2^\kappa$ we have our claim.
	
	We will proceed by indction on $\alpha < 2^\kappa$. Before starting the induction we enumerate
	$$\mathcal{P}(\omega) = \{M_\alpha \colon \alpha < 2^\kappa, \alpha \equiv 0 (mod 2)\}$$
	and enumerate the family of all partitions of $\kappa$ 
	$$\{\{D^\alpha_t \colon t  in [\kappa]^{<\omega}\} \colon \alpha < 2^\kappa, \alpha \equiv 1 (mod 2)\}$$
	in such a way that each partition is listed $2^\kappa$ times.
	
	We will define $\mathcal{F}^\alpha_{\mu, t}$, $\mathcal{F}^\alpha_{\xi, t}$ and $I_\alpha$ as follows. For $\alpha = 0$ take
	$$\mathcal{F}^0_{\mu, t} = [\mathcal{FR}(\kappa), \{\bigcup_{s \supseteq \hat{\varphi}(p)} A^0_\xi\cap E^0_s \colon \xi< \mu, p \in [\kappa^+]^{< \omega}, \xi \in p\}]$$
	$$\mathcal{F}^0_{\xi, t} = [\mathcal{FR}(\kappa), \{\bigcup_{s \supseteq \hat{\varphi}(p)} A^0_\eta\cap E^0_s \colon \eta< \xi, p \in [\kappa^+]^{< \omega}, \eta \in p\}, \{A^0_\xi \cap E^0_t\},$$ $$ \{A^0_\zeta\cap E^0_s \colon s \not\supseteq \hat{\varphi}, \xi < \zeta < \mu, p \in [\kappa^+]^{< \omega}, \xi \in p\}]$$
	and 
	$$I_0 = 2^\kappa \setminus \mu.$$
	Further
	\begin{itemize}
		\item [(i)] If $\beta < \alpha$ then $\mathcal{F}^\beta_{\xi, t} \subseteq \mathcal{F}^\alpha_{\xi, t}$ and $I_\beta \supseteq I_\alpha$ for any $\xi \leqslant \mu$ and $t \in [\kappa]^{<\omega}$,
		\item [(ii)] $\mathcal{F}^\gamma_{\xi, t} = \bigcup_{\beta < \gamma} \mathcal{F}^\beta_{\xi, t}$ and $I_\gamma = \bigcap_{\beta < \gamma} I_\beta$ for limit $\gamma$,
		\item [(iii)] $|2^\kappa\setminus I_\alpha| \leqslant |\alpha| + \mu$ and $I_\alpha \setminus I_{\alpha+1}$ is finite,
		\item [(iv)] each  $\mathcal{F}^\alpha_{\xi, t}$, $\xi \leqslant \mu, t \in [\kappa]^{<\omega}$ isa filter on $\kappa$ and the matrix 
			$$\{E_t^i \colon t \in [\kappa]^{<\omega}, i \in I_\alpha\} \cup \{A_\eta^i \colon \eta < \kappa^+, i \in I_\alpha\}$$
			of remaining step families is independent with respect to  $\mathcal{F}^\alpha_{\xi, t}, \hat{\varphi}$,
		\item [(v)]	For any $M \in  \mathcal{F}^\alpha_{\xi, t}$and $\zeta < \xi \leqslant \mu$ the set $\{s \in [\kappa]^{< \omega} \colon M \in  \mathcal{F}^\alpha_{\xi, s}\}$ has cardinality $\kappa$,
		\item [(vi)] There exists $M \in \mathcal{F}^\alpha_{\xi, s}$ such that $\kappa \setminus M \in \mathcal{F}^\alpha_{\xi, t}$ for any $t \in [\kappa]^{<\omega}$ whenever $\zeta < \xi \leqslant \mu$ and $s \in [\kappa]^{<\omega}$,
		\item [(vii)] If $\alpha \equiv 0 (mod 2)$ then either $M_\alpha \in \mathcal{F}^{\alpha+1}_{\xi, s}$ or $\kappa \setminus M_\alpha \in \mathcal{F}^{\alpha+1}_{\xi, s}$ for any $\xi \leqslant \mu$ and $t \in [\kappa]^{<\omega}$,
		\item [(viii)] If $\alpha \equiv 1 (mod 2)$and if 
		$\{D^\alpha_t \colon t \in [\kappa]^{<\omega}\}$ is such that for each $t \in [\kappa]^{< \omega}$ we have $\kappa \setminus D^\alpha_t \in \mathcal{F}^\alpha_{\xi, s}$ and for each $\zeta < \mu$ the set $\{s \in [\kappa]^{<\omega}\colon D^\alpha_t \in \mathcal{F}^\alpha_{\xi, s}\}$ has cardinality $\kappa$ then for each $\kappa-$discrete set $Y$ if
		$$Y \subseteq \bigcap_{\xi < \mu} \overline{\{\mathcal{F}^\alpha_{\xi, t} \colon t \in [\kappa]^{<\omega}\}} \setminus \mathcal{F}^\alpha_{\mu, t}$$
		then $\mathcal{F}^\alpha_{\mu, t} \not \in \overline{Y}$.
 		\end{itemize}
 	Notice that $(v)$ gives us that $\mathcal{F}^\alpha_{\xi, t} \in \overline{\{\mathcal{F}^\alpha_{\xi, t} \colon s \in [\kappa]^{<\omega}\}}$ for any $\xi < \mu$ and $ t \in [\kappa]^{<\omega}$ while $(vi)$ gives us that $\mathcal{F}^\alpha_{\zeta, s} \not \in \overline{\{\mathcal{F}^\alpha_{\xi, t} \colon t \in [\kappa]^{< \omega}\}}$ whenever $\zeta < \xi < \mu$ and $s \in [\kappa]^{<\omega}$.
 	Further, $(VII)$ assures that $\mathcal{G}_{\xi, t}$ will be an ultrafilter for any $\xi \leqslant \mu$ and $t \in [\kappa]^{<\omega}$ and $(viii)$ assures that $(3)$ will be fulfilled.
 	
 	As usual in these constructions, there is no problem at limits. Thus, we proceed to describe the successor steps.
 	
 	If $\alpha \equiv 0 (mod 2)$. If $\mathcal{R} = [\mathcal{F}^\alpha_{\xi, t}, \{M_|alpha\}]$ is a proper filter and the matrix 
 		$$\{E_t^i \colon t \in [\kappa]^{<\omega}, i \in I_\alpha\} \cup \{A_\eta^i \colon \eta < \kappa^+, i \in I_\alpha\}$$
 		is indepenednt with respect to $\mathcal{R}, \hat(\varphi)$ then we set $$\mathcal{F}^{\alpha+1}_{\xi, t} = [\mathcal{F}^\alpha_{\xi, t}, \{M_\alpha\}]$$
 		and $I_{\alpha+1} = I_\alpha$.
 
 		Otherwise, fix $n \in \omega$, distinct $i_k \in I_\alpha$ and $\hat{\varphi}(p_k)\subseteq t_k,$ for $ k < n$, such that $$\kappa \setminus [Z_\alpha \cap \bigcap_{k=0}^{n-1}(\bigcap_{\eta\in p_k}A^{i_k}_{\eta} \cap E^{i_k}_{t_k})]\in \mathcal{F}^\alpha_{\xi, \zeta}.$$ Then, put $$\mathcal{F}^{\alpha}_{\xi, \zeta} = [\mathcal{F}^{\beta}_{\xi, \zeta}, \{A^{i_k}_{p_k} \colon 0 \leqslant k \leqslant n-1\}, \{E^{i_k}_{t_k} \colon 0 \leqslant k \leqslant n-1\}]$$ $$I_{\alpha} = I_{\beta}\setminus \{i_k \colon 0 \leqslant k \leqslant n-1\}.$$ Then $\kappa \setminus Z_\alpha \in \mathcal{F}^{\alpha}_{\xi, \zeta}.$ We leave it to the readre to chceck that the condition $(iv)$ holds.
 		
 		If $\alpha \equiv 1 (mod 2)$. Assume that for each $t \in [\kappa]^{<\omega}$ we have $\kappa \setminus D^\alpha_t \in \mathcal{F}^\alpha_{\mu, s}, s \in [\kappa]^{<\omega}$ and for each $\zeta < \mu$ the set
 		$$\{s \in [\kappa]^{<\omega} \colon D^\delta_{t} \in \mathcal{F}^\alpha_{\zeta, s}\}$$
 		has cardinality $\kappa$.
 		
 		Otherwise, we set $\mathcal{F}^{\alpha+1}_{\xi, t} = \mathcal{F}^\alpha_{\xi, t}$ and $I_{\alpha+1} = I_\alpha$.
 		Choose $i \in I_\alpha$ and set $I_{\alpha+1} = I_\alpha\setminus \{i\}$.
 		Since $D^\alpha_t \in \mathcal{F}^\alpha_{\zeta, s}$ for $\kappa$ many $s \in [\kappa]^{<\omega}$ the intersection
 		$$D^\alpha_t \cap (A^i_\zeta \cap E^i_s)$$
 		has cardianlity $\kappa$. Thus, we put
 		 $$\mathcal{F}^{\alpha+1}_{\mu, t} = [\mathcal{F}^\alpha_{\mu, t}, \{\bigcup_{s \in [\kappa]^{<\omega}}(\kappa\setminus D^\alpha_t)\cap(A^i_\zeta \cap E^i_s) \colon s \subseteq t\}].$$
 		 To define $\mathcal{F}^{\alpha+1}_{\xi, t}$ for $\xi < \mu$ we pick one $i \in I_\alpha$, put $I_{\alpha+1} = I_\alpha \setminus \{i\}$ and put
 		 $$\mathcal{F}^{\alpha+1}_{\xi, t} = [\mathcal{F}^\alpha_{\xi, t}, \{\bigcup_{s \in [\kappa]^{<\omega}}(D^\alpha_t)\cap(A^i_\zeta \cap E^i_s) \colon s \subseteq t\}].$$
 		 To see that $(viii)$ is fulfilled take a $\kappa-$discrete set $$Y \subseteq \bigcap_{\xi < \mu}\overline{\{\mathcal{F}^{\alpha+1}_{\xi, t} \colon t \in [\kappa]^{<\omega}\}}\setminus \mathcal{F}^{\alpha+1}_{\mu, t}.$$
 		 Since $Y \subseteq \bigcap_{\xi < \mu}\overline{\{\mathcal{F}^{\alpha+1}_{\xi, t} \colon t \in [\kappa]^{<\omega}\}}$  the set $\{s \in [\kappa]^{<\omega} \colon D^\alpha_t \in \mathcal{F}^\alpha_{\xi, t}\}$ has cardinality $\kappa$. Hence the assumptions in $(viii)$ are fulfilled. Then $\kappa\setminus D^\alpha_t)\cap(A^i_\zeta \cap E^i_s) \in \mathcal{F}^{\alpha+1}_{\xi, t}$ for all $t \in [\kappa]^{<\omega}$ and $s \subseteq t$. Hence $\mathcal{F}^{\alpha+1}_{\mu, t} \not in \overline{Y}$.
 		 
 		 To see that $(iv)$ holds notice that each element of $\mathcal{F}^\alpha_{\mu, t}$ is of the form 
 		 $$B \cap ((\kappa\setminus D^\alpha_t)\cap (A^i_\xi \cap E^i_s))$$
 		 for sme $B \in \mathcal{F}^\alpha_{\mu, t}$. But the assumption $\kappa \setminus D^\alpha_t \in \mathcal{F}^\alpha_{\mu, t}$ and by $(a)$ and $(d)$ $\{E^i_s \colon s \in [\kappa]<{<\omega}\}$ is a partition of $\kappa$. Hence $A^i_\xi \cap E^i_s \in \mathcal{F}^\alpha_{\mu, t}$
 for $s \subseteq t$.
 Hence the condition $(iv)$	for $\mathcal{F}^{\alpha+1}_{\mu, t}$ follows from the condition $(iv)$ for $\mathcal{F}^\alpha_{\mu, t}$.
 In order tp show that $(iv)$ is fulfilled for $\mathcal{F}^\alpha_{\xi, t}$ is similar to the above argumentation and it is left to the readre. The proof is complete.
	\end{proof}
	\\

 \noindent
\textbf{Open problem} Does there exist a chain in the Rudin-Frol\'ik order of $\beta\kappa\setminus \kappa$ of length $\mu$, $\kappa \leqslant \mu \leqslant 2^\kappa$ for singular $\kappa> \omega$ without a lower bound?

{\sc Joanna Jureczko}
\\
Wroc\l{}aw University of Science and Technology, Wroc\l{}aw, Poland
\\
{\sl e-mail: joanna.jureczko@pwr.edu.pl}

\end{document}